
\input amstex
\documentstyle{amsppt}

\magnification=\magstep1
\hsize=6.5truein
\vsize=9truein

\font \smallrm=cmr10 at 10truept
 at 7truept
\font \smallbf=cmbx10 at 10truept
 at 10truept
\font \smallsl=cmsl10 at 10truept

\baselineskip=.15truein

\def \N {\Bbb N}

\def \Id {\mathop{\hbox{\rm Id}}\nolimits}
\def \End {\mathop{\hbox{\rm End}}\nolimits}\def \mod {\mathop{\hbox{\rm
mod}}\nolimits}

\def \e {\frak e}
\def \g {\frak g}
\def \kh {k[[h]]}

\def \ug {U(\g)}
\def \uhg {U_h(\g)}

\document

\topmatter

{\ }

\vskip-57pt

\hfill   {{\smallsl Journal of Pure and Applied Algebra\/}
{\smallbf 161}  {\smallrm (2001), 295--307}}
\hskip19pt   {\ }
                                      \par
\hfill   {\smallrm DOI: 10.1016/S0022-4049(00)00099-2}
\hskip19pt   {\ }

\vskip31pt

\title
   Tressages des groupes de Poisson formels \`a dual quasitriangulaire
\endtitle

\author
        Fabio Gavarini${}^\dag$ \, ,  \;  Gilles Halbout${}^\ddag$
\endauthor

\leftheadtext{ Fabio Gavarini, Gilles Halbout }
\rightheadtext{ Tressages des groupes de Poisson formels \`a dual
quasitriangulaire }

\affil
 ${}^\dag$ \! Universit\`a di  \hbox{ Roma ``Tor Vergata'',
Dipartimento di Matematica }  --  Roma, Italie   \\
  ${}^\ddag$ \! Institut de  \hbox{ Recherche Math\'ematique Avanc\'ee,
ULP--CNRS   --   Strasbourg, } France  \\
\endaffil

\address\hskip-\parindent
  ${}^\dag$ \! Universit\`a degli Studi di Roma ``Tor Vergata''   ---
Dipartimento di Matematica   \newline
  Via della Ricerca Scientifica, 1   ---   I-00133 Roma, Italie   ---
e-mail: gavarini\@{}mat.uniroma2.it  \newline
     \newline
  ${}^\ddag$ \! Institut de Recherche Math\'ematique Avanc\'ee, Universit\'e
  Louis Pasteur
  et C.N.R.S.   ---
  \newline
  7, rue Ren\'e{} Descartes   ---   67084 STRASBOURG Cedex, France ---
  e-mail:  halbout\@math.u-strasbg.fr
\endaddress

\abstract
   In [D3], Drinfeld constructs a Quantum Formal Series Hopf Algebra
(QFSHA) $U'_h$  starting from a Quantum Universal Enveloping Algebra
(QUEA) $U_h$.  In this paper, we prove that if  $ \big( U_h, R \big) $
is any quasitriangular QUEA, then  $ \Big( {U'_h}, \hbox{Ad}(R)
{\big\vert}_{{U'_h} \otimes {U'_h}} \Big) $  is a braided QFSHA.
As a consequence, we prove that if $ \g $  is a quasitriangular
Lie bialgebra over a field  $ k $  of characteristic zero and
$\g^*$ is its dual Lie bialgebra, the algebra of functions $ F[[\g^*]] $ on
the formal group associated to $\g^*$ is a
braided Hopf algebra.  This result is a consequence of
the existence of a quasitriangular quantization  $ (U_h,R) $
of  $ U(\g) $  and of the fact that  $ {U'_h} $  is a
quantization of  $ F[[\g^*]] $.
\endabstract

\endtopmatter

  \footnote""{ ${}^\dag$ \! Le premier auteur a \'et\'e en partie
financ\'e par une bourse du  {\it Consiglio Nazionale delle Ricerche\/}
\, (Italie) }

\vskip11pt

\centerline {\bf  Introduction }

\vskip10pt

   Soit  $ \g $  une big\`ebre de Lie sur un corps  $ k $
de caract\'eristique z\'ero et $ \g^* $  sa big\`ebre de
Lie duale (topologique, en g\'en\'eral).  L'alg\`ebre
$ F[[\g^*]] $  des fonctions sur
le groupe de Poisson formel associ\'e \`a  $ \g^* \, $
est une  $ k $-alg\`ebre  de Hopf topologique.
Dans ce travail, nous d\'emontrons
que la donn\'ee d'une structure
quasitriangulaire sur  $ \g $  (c'est-\`a-dire
d'une $ r $-matrice,  $ r \in \wedge^2 \g $,  solution
de l'\'equation de Yang-Baxter
classique) induit un tressage sur l'alg\`ebre de Hopf $ F[[\g^*]] $
(la d\'efinition d'un tressage est donn\'ee au paragraphe $1$).
Nous \'etendons ainsi au cas g\'en\'eral les r\'esultats
obtenus pour les alg\`ebres de Kac-Moody de type fini ou affine
par Reshetikhin [Re] et le premier auteur [G1,G2].

\smallskip

Notre d\'emonstration repose sur les quantifications
des alg\`ebres enveloppantes et
des alg\`ebres de fonctions sur les groupes formels.
L'alg\`ebre enveloppante $\ug$ d'une big\`ebre de Lie
$ \g $  est une alg\`ebre de Hopf-co-Poisson.
Dans ce cadre,
Etingof et Kazhdan [EK] ont montr\'e l'existence d'une
quantification de $\ug$, c'est-\`a-dire
d'une  $ k[[h]] $-alg\`ebre  de Hopf topologique $\uhg$ v\'erifiant:
\itemitem{(a)} les $k[[h]]$-modules $\uhg$ et $\ug[[h]]$ sont isomorphes;
\itemitem{(b)} l'alg\`ebre de Hopf-co-Poisson $\uhg \otimes_{k[[h]]} k$
obtenue par sp\'ecialisation \`a $h=0$ est isomorphe \`a $\ug$.

A partir d'une big\`ebre de Lie $\g$, on peut construire {\it a priori}
plusieurs quantifications $\uhg$. A l'int\'erieur
de chacune d'elles, Drinfeld [D3] construit une sous-alg\`ebre de Hopf
$F_h[[\g^*]]$ dont la sp\'ecialisation \`a $h=0$ est isomorphe \`a
l'alg\`ebre de Hopf-Poisson $F[[\g^*]]$.

\smallskip

Supposons maintenant que $\g$ soit en outre quasitriangulaire, munie
d'une  $ r $-matrice  $ r \in \wedge^2\g \subseteq \g \otimes \g $.
Etingof et Kazhdan ont montr\'e, dans [EK], que cette structure
qua-sitriangulaire pouvait elle aussi \^etre quantifi\'ee:
une des quantifications  $ \uhg $  poss\`ede une  $ R $-matrice
universelle  $ R_h \in \uhg \otimes \uhg $  (produit tensoriel
topologique) qui, dans l'identification de  $ k[[h]] $-modules
$ \uhg \otimes \uhg \simeq \left( \ug \otimes \ug\right)[[h]] $,
s'\'ecrit sous la forme  $ R_h\equiv 1 \otimes 1 + h r~ (\mod h^2) $.
Nous prouvons dans ce cadre que l'action de  $ R_h $  par automorphisme
int\'erieur dans  $ \uhg \otimes \uhg $  stabilise la sous-alg\`ebre
$ F_h[[\g^*]] \otimes F_h[[\g^*]] $  et induit par sp\'ecialisation
un op\'erateur  $ {\frak R}_0 $  sur  $ F[[\g^*]] \otimes F[[\g^*]] $.
Les propri\'et\'es alg\'ebriques de la  $ R $-matrice  universelle
$ R_h $  font que ${\frak R}_0$ est un op\'erateur de tressage, ce
qui d\'emontre le r\'esultat.

\smallskip

   Dans la premi\`ere partie de ce papier, nous rappellerons les
d\'efinitions et notions utiles pour notre travail.  Dans la seconde
partie, nous \'enoncerons pr\'ecis\'ement les r\'esultats principaux
et donnons le sch\'ema de leurs d\'emonstrations.  Dans la troisi\`eme
partie, on trouvera la preuve, essentiellement combinatoire, du
th\'eor\`eme technique 2.1.

\vskip10pt

\centerline { REMERCIEMENTS }

\vskip7pt

 Les auteurs tiennent \`a remercier M.~Rosso et
C.~Kassel pour de nombreux entretiens.

\vskip1,6truecm

\centerline {\bf \S\; 1. \ D\'efinitions et rappels }

\vskip10pt

  {\bf 1.1  Les objets classiques.}  Fixons un corps  $ k $  de
caract\'eristique z\'ero qui sera le corps de base de tous les
objets classiques (alg\`ebres et big\`ebres de Lie, alg\`ebres
de Hopf, etc.) que nous introduirons.
                                                  \par
   Suivant [D1], nous appelons big\`ebre de Lie une paire
$ \, (\g, \delta_\g) \, $  o\`u  $ \g $  est une alg\`ebre de Lie et
$ \, \delta_\g \colon \, \g \rightarrow \g \otimes \g \, $
est une application lin\'eaire antisym\'etrique   --- dite cocrochet
de Lie ---   telle que son dual  $ \, \delta_\g^* \colon \, \g^*
\otimes \g^* \rightarrow \g^* \, $  soit un crochet de Lie et que
$ \delta_g $  elle-m\^eme soit un 1-cocycle de  $ \g $  \`a valeurs
dans  $ \g \otimes \g $.  Le dual lin\'eaire  $ \g^* $ de  $ \g $
est alors \`a son tour une big\`ebre de Lie (topologique
par rapport \`a la topologie faible lorsque
dim$ (\g) = \infty \, $).  Suivant [D3], \S 4, nous appelons
big\`ebre de Lie quasitriangulaire un couple  $ \, (\g,r) \, $
v\'erifiant les propri\'et\'es suivantes:  $ \, r \in \g \otimes \g \, $
est solution de l'\'equation de Yang-Baxter classique (CYBE) dans
$ \, \g \otimes \g \otimes \g \, $:
$$ \, [r_{1{}2},r_{1{}3}] +
[r_{1{}2}, r_{2{}3}] + [r_{1{}3},r_{2{}3}] = 0 \, $$
et  $ \g $
est une big\`ebre de Lie munie du cocrochet  $ \, \delta =
\delta_\g \, $  d\'efini par  $ \, \delta(x) = [x\otimes 1 +
1 \otimes x,r] \, $;  l'\'el\'ement  $ r $  est alors appel\'e
la  $ r $-matrice  de  $ \g $.
                                                  \par
   Si  $ \g $  est une alg\`ebre de Lie, son alg\`ebre
enveloppante universelle  $ U(\g) $  est une alg\`ebre de Hopf;
si de plus  $ \g $  est une big\`ebre de Lie, alors  $ U(\g) $
est en fait une alg\`ebre de Hopf-co-Poisson ([D2]).
                                                  \par
   Soit  $ \g $  une alg\`ebre de Lie quelconque. Comme  $ U(\g ) $
est une alg\`ebre de Hopf, son dual est une alg\`ebre de Hopf formelle
([Di], Ch.~1); on appelle alors alg\`ebre des fonctions sur le groupe
formel associ\'e \`a  $ \g $,  ou tout simplement groupe formel associ\'e
\`a  $ \g $,  cette alg\`ebre de Hopf formelle  $ \, F[[\g]] = \!\!
{\phantom{\big(} U(\g) \phantom{\big)}}^{\!\!\!*} \, $.  Si  $ G $
est un groupe alg\'ebrique affine connexe d'alg\`ebre de Lie
$ \g $, si $ F[G] $  est l'alg\`ebre de Hopf des fonctions
r\'eguli\`eres sur  $ G $, et si  $ {\frak m}_e $  est l'id\'eal
maximal dans  $ F[G] $  des fonctions qui s'annulent au point unit\'e
$ \, e \in G \, $,  alors l'alg\`ebre de Hopf formelle  $ \, F[[\g]] \, $
s'identifie \`a la compl\'etion  $ {\frak m}_e $-adique  de  $ F[G] $
([On], Ch.~I).  Lorsque, de plus,  $ \g $  est une big\`ebre de Lie,
$ F[[\g]] $  est en fait une alg\`ebre de Hopf-Poisson formelle
([CP], \S 6.2.A).

\vskip7pt

  {\bf 1.2  Tressages et quasitriangularit\'e.}  Soit  $ H $  une
alg\`ebre de Hopf dans une cat\'egorie tensorielle  $ ({\Cal A},\otimes) $
([CP], \S 5):  $ H $  est dite tress\'ee ([Re], {Def.} 2)
s'il existe un automorphisme  ${\frak R} $  de l'alg\`ebre
$ H \otimes H $, appel\'e op\'erateur de tressage de
$ H $,  diff\'erent de  la volte $ \; \sigma \colon a \otimes b \mapsto
b \otimes a \; $,  et v\'erifiant
  $$  \displaylines{
   {\frak R} \circ \Delta = \Delta^{\hbox{\smallrm op}},  \cr
   (\Delta \otimes {\Id}) \circ {\frak R} = {\frak R}_{13} \circ
{\frak R}_{23} \circ (\Delta \otimes {\Id}) \; ,  \qquad
({\Id} \otimes \Delta) \circ {\frak R} = {\frak R}_{13} \circ
{\frak R}_{12} \circ ({\Id} \otimes \Delta)  \cr }  $$
o\`u   $ \Delta^{\hbox{\smallrm op}} = \sigma \circ \Delta $  et
$ {\frak R}_{12} $,  $ {\frak R}_{13} $  et  $ {\frak R}_{23} $  sont
les automorphismes de  $ H \otimes H \otimes H $  d\'efinis par  $ \,
{\frak R}_{12} = {\frak R} \otimes {\Id} \, $,  $ \,
{\frak R}_{23} = {\Id} \otimes {\frak R} \, $,  $ \,
{\frak R}_{13} = (\sigma \otimes {\Id}) \circ ({\Id} \otimes {\frak R})
\circ (\sigma \otimes  {\Id}) \, $.
                                                  \par
   Si l'alg\`ebre $ H $  est une alg\`ebre de Hopf-Poisson,
   nous dirons que $H$ est tress\'ee  en tant
qu'alg\`ebre de Hopf-Poisson  si elle est tress\'ee  en tant qu'alg\`ebre
de Hopf et si son tressage est un automorphisme
d'alg\`ebre de Poisson.

\vskip3pt

   Si la paire  $ (H, {\frak R}) $  est une alg\`ebre tress\'ee, il
r\'esulte de la d\'efinition que  $ {\frak R} $  v\'erifie l'\'equation
de Yang-Baxter quantique (QYBE) dans
$ \End(H^{\otimes 3}) $:
  $$  {\frak R}_{12} \circ {\frak R}_{13} \circ {\frak R}_{23}
= {\frak R}_{23} \circ {\frak R}_{13} \circ {\frak R}_{12} . $$
Alors, pour tout  $ \, n \in \N \, $, le groupe des
tresses  $ {\Cal B}_n $  agit sur  $ H^{\otimes n} $, ce qui
permet de construire des invariants de n\oe{}uds ([Tu]).

\vskip3pt

  Une alg\`ebre de Hopf  $ H $   dans une cat\'egorie tensorielle  est
dite quasitriangulaire ([D3], [CP]) s'il existe un \'el\'ement
inversible  $ \, R \in H \otimes H \, $,  appel\'e  $ R $-matrice  de
$ H $,  tel que
  $$  \displaylines{
  {\hbox{\rm Ad}}(R) (\Delta
(a)) =  R \cdot \Delta (a) \cdot R^{-1} = \Delta^{\hbox{\smallrm op}}(a),  \cr
   (\Delta \otimes {\Id}) (R) = R_{13} R_{23} \; ,  \qquad
({\Id} \otimes \Delta) (R) = R_{13} R_{12}  \cr }  $$
o\`u  $ \, R_{12}$, $R_{13}$ et $ R_{23}$ sont les
\'el\'ements de  $H^{\otimes 3} $ d\'efinis par  $ \, R_{12}
= R \otimes 1 \, $,  $ \, R_{23} = 1 \otimes R \, $ et  $ \, R_{13} =
(\sigma \otimes {\Id}) (R_{23})\, $.
Il r\'esulte classiquement des identit\'es ci-dessus que  $ R $
v\'erifie la QYBE dans  $ H^{\otimes 3} $,  c'est-\`a-dire,
  $$  R_{12} R_{13} R_{23} = R_{23} R_{13} R_{12}  .  $$
Les produits tensoriels de  $ H $-modules  sont alors munis d'une
action du groupe des tresses.  En outre, il est clair que si
$ (H,R) $  est quasitriangulaire, alors  $ \big( H, {\hbox{\rm
Ad}}(R) \big) $  est tress\'ee.

\vskip7pt

  {\bf 1.3  Les objets quantiques.}  Soit  $ {\Cal A} $  la cat\'egorie
dont les objets sont les  $ \kh $-modules  topologiquement libres et
complets au sens  $ h $-adique, et les morphismes sont les
applications  $ \kh $-lin\'eaires continues. Pour tous  $ V $,  $ W $
dans  $ {\Cal A} $,  d\'efinissons  $ \, V \otimes W \, $ comme
\'etant  la limite projective des  $ \kh \big/ (h^n) $-modules
$ \, \big( V / h^n V \big) \otimes_{\kh / (h^n)} \big( W / h^n W \big)
\, $:  on munit ainsi  $ {\Cal A} $  d'une structure de cat\'egorie
tensorielle; en particulier pour
tous espaces vectoriels  $ V_0 $ et  $ W_0 $  sur  $ k $  on a  $ \,
V_0[[h]] \otimes W_0[[h]] \cong (V_0 \otimes W_0)[[h]] \, $,  o\`u
$ V_0 \otimes W_0 $  est un produit tensoriel topologique si
$ V_0 $ et  $ W_0 $  sont des espaces vectoriels topologiques ([EK],
\S 7).  Reprenant la d\'efinition de  Drinfeld ([D3]), on
appelle alg\`ebre enveloppante universelle quantifi\'ee (QUEA) toute
alg\`ebre de Hopf dans la cat\'egorie  $ {\Cal A} $ dont la limite
semi-classique, c'est-\`a-dire la sp\'ecialisation en  $ \, h = 0 $,
est l'alg\`ebre enveloppante universelle d'une big\`ebre de Lie.
En particulier, toute quantification d'une big\`ebre de Lie est
une QUEA.  De m\^eme, on appelle alg\`ebre de Hopf des s\'eries
formelles quantiques (QFSHA), toute alg\`ebre de Hopf dont la
limite semi-classique est l'alg\`ebre des fonctions sur un
groupe formel.

\vskip3pt

   Dans la suite, nous aurons besoin du r\'esultat suivant:

\vskip7pt

\proclaim{Th\'eor\`eme 1.4}  ([EK])  Toute big\`ebre
de Lie $ \g $  admet une quantification $ \, \uhg $.
Si  $( \g, r) $  est de plus quasitriangulaire, alors
il existe une quantification  $ \, \uhg \, $ et un
\'el\'ement  $ \, R_h \in U_h(\g) \otimes U_h(\g) \, $
tels que  $ \, \big( U_h(\g), R_h \big) \, $  soit une
alg\`ebre de Hopf quasitriangulaire et  $ \, R_h \in 1
\otimes 1  + r \, h +  h^2 \left( \ug \otimes \ug\right)[[h]] $.
\endproclaim

\vskip7pt

  {\bf 1.5  Le foncteur de Drinfeld.}  Soit  $ H $  une alg\`ebre de
Hopf sur  $ \kh $.  Pour tout  $ \, n \in \N $,  on d\'efinit
$ \; \Delta^n \colon \, H \longrightarrow H^{\otimes n} \; $  par  $ \,
\Delta^0 = \epsilon $,  $ \, \Delta^1 = {\Id}_{\scriptscriptstyle
H} $  et  $ \, \Delta^n = \big( \Delta \otimes
{\Id}_{\scriptscriptstyle H}^{\otimes (n-2)} \big) \circ
\Delta^{n-1} \, $  si  $ \, n > 2 $.  Pour tout sous-ensemble
ordonn\'e  $ \, \Sigma = \{i_1, \dots, i_k\} \subseteq \{1, \dots,
n\} \, $  avec  $ \, i_1 < \dots < i_k \, $,  \, on d\'efinit
l'homomorphisme  $ \; j_{\scriptscriptstyle \Sigma} \colon \,
H^{\otimes k} \longrightarrow H^{\otimes n} \; $  par  $ \;
j_{\scriptscriptstyle \Sigma} (a_1 \otimes \cdots \otimes a_k) =
b_1 \otimes \cdots \otimes b_n \; $  avec  $ \, b_i = 1 \, $  si
$ \, i \notin \Sigma \, $  et  $ \, b_{i_m} = a_m \, $  pour
$ \, 1 \leq m \leq k \, $;  on pose alors  $ \; \Delta_\Sigma =
j_{\scriptscriptstyle \Sigma} \circ \Delta^k \, $.  On d\'efinit
aussi  $ \; \delta_n \colon \, H \longrightarrow H^{\otimes n}
\; $  par  $ \; \delta_n = \sum_{\Sigma \subseteq \{1, \dots, n\}}
{(-1)}^{n-|\Sigma|} \Delta_\Sigma \, $,  pour tout  $ \, n \in \N_+
\, $,  et plus  g\'en\'eralement, pour tout  $ \, \Sigma = \{i_1, \dots,
i_k\} \subseteq \{1, \dots, n\} \, $,  avec  $ \, i_1 < \cdots < i_k \, $,
on pose
  $$  \delta_\Sigma = {\textstyle \sum_{\Sigma' \subseteq \Sigma}}
{(-1)}^{|\Sigma|-|\Sigma'|} \, \Delta_{\Sigma'}
=j_\Sigma \circ \delta_k .   \eqno (1.1)  $$
En particulier,  $ \, \delta_{\{1, \dots, n\}} = \delta_n \, $.
Gr\^ace au principe d'inclusion-exclusion, ceci \'equivaut \`a
  $$  \Delta_\Sigma = {\textstyle \sum_{\Sigma' \subseteq \Sigma }}
\delta_{\Sigma'}   \eqno (1.2)  $$
pour tout  $ \, \Sigma = \{i_1, \dots, i_k\} \subseteq \{1, \dots, n\}
\, $  avec  $ \, i_1 < \cdots < i_k \, $.  Enfin on d\'efinit le
sous-espace de $H$
  $$  H' = \big\{\, a \in H \,\big\vert\, \delta_n(a) \in h^n
H^{\otimes n} \, \big\} \, ,  $$
que nous consid\'ererons muni de la topologie
induite.  Nous avons alors:

\vskip7pt

\proclaim{Th\'eor\`eme 1.6}  ([D3])  Si  $ H $  est
une QUEA, alors  $ \, H' \, $  est une QFSHA.
\hfill\break
   \indent   Si, de plus, la limite semi-classique de $H$ est
l'alg\`ebre de Hopf-co-Poisson $\ug$, alors
la limite semi-classique de
$ {\phantom{\big(} U_h(\g) \phantom{\big)}}^{\!\!\prime} $  est
l'alg\`ebre de Hopf-Poisson toplogique $ \, F[[\g^*]] \, $.
\endproclaim

\vskip2,1truecm

\centerline {\bf \S \; 2. \  Les r\'esultats principaux }

\vskip10pt

Les r\'esultats de cet article sont des cons\'equences
du th\'eor\`eme suivant dont
la preuve sera donn\'ee dans le paragraphe 3.
%
%
 \eject

\proclaim{Th\'eor\`eme 2.1}  Soit  $ H $  une alg\`ebre de Hopf
quasitriangulaire dans la cat\'egorie  $ {\Cal A} $,  et soit  $ R $  sa
$ R $-matrice.  Alors l'automorphisme interieur  $ \, {\hbox{\rm Ad}}
(R) \colon \, H \otimes H \rightarrow H \otimes H \, $  se
restreint en un automorphisme de  $ \, H' \otimes H' $ et la paire  $ \,
\Big( H', \, {\hbox{\rm Ad}}(R) {\big\vert}_{H' \otimes H'} \Big) \, $
est ainsi une alg\`ebre de Hopf tress\'ee.
\endproclaim

On d\'eduit de ce th\'eor\`eme
une interpr\'etation g\'eom\'etrique de la  $ r $-matrice  classique:

\vskip7pt

\proclaim{Th\'eor\`eme 2.2}  Pour toute  big\`ebre de Lie
quasitriangulaire $ \g $, l'alg\`ebre de Hopf-Poisson topologique
$ \, F[[\g^*]] \, $  est tress\'ee.
Plus pr\'ecis\'ement, on peut trouver
une alg\`ebre de Hopf tress\'ee qui quantifie
l'alg\`ebre $ F[[\g^*]]$
et dont l'op\'erateur de tressage se sp\'ecialise en celui
de  $ F[[\g^*]] $.
\endproclaim

\demo{Preuve}  Soit  $ r $  la  $ r $-matrice  de  $ \g $.  D'apr\`es
le th\'eor\`eme 1.4, il existe une QUEA quasitriangulaire  $ \, \big(
\uhg, R_h \big) \, $  dont la limite semi-classique est
$ \, \big( \ug, \, r \big) \, $. D'apr\`es le th\'eor\`eme 1.6, la limite
semi-classique de  $ \uhg' $  est  $ F[[\g^*]] $.  Soit  $ \, {\frak R}_h
= {\hbox{\rm Ad}}(R_h)$,
l'automorphisme int\'erieur d\'efini par
$R_h$. Le th\'eor\`eme 2.1 nous assure que
 $ \, \Big( \!\! {\phantom{\big(} \uhg \phantom{\big)}}^{\!\!\prime},
\, {\frak R}_h{\big\vert}_{{\uhg}' \otimes {\uhg}'} \, \Big) \, $  est
une alg\`ebre de Hopf tress\'ee. Sa limite semi-classique  $ \,
\bigg( F[[\g^*]], \, \left( {\frak R}_h{\big\vert}_{{\uhg}' \otimes
{\uhg}'} \right) {\Big\vert}_{h=0} \bigg) \, $  est donc tress\'ee elle-aussi.
\smallskip
Enfin, le crochet de Poisson de  $ F[[\g^*]] $
est donn\'e par  $ \, \{ a, b \} = {\big( [ \alpha, \beta ] \big/ h
\big)}{\big\vert}_{h=0} \, $  pour tout  $ \, a $,  $ b \in F[[\g^*]] \, $
et  $ \, \alpha $,  $ \beta \in \!\!\! {\phantom{\big(} \uhg
\phantom{\big)}}^{\!\!\prime} \, $  tels que  $ \, \alpha{\vert}_{h=0} = a
\, $,  $ \, \beta{\vert}_{h=0} = b $.
Ainsi, puisque $ {\frak R}_h $  est un automorphisme d'alg\`ebre,
sa restriction $ \left( {\frak
R}_h{\big\vert}_{{\uhg}' \otimes {\uhg}'} \right) {\Big\vert}_{h=0} $  est
aussi un automorphisme d'alg\`ebre de Poisson.   $ \square $
\enddemo

\vskip7pt

   Le th\'eor\`eme ci-dessus a une autre cons\'equence:
soient  $ \g $  et  $ \g^* $  comme ci-dessus, soit  $ {\frak R} $  le
tressage de  $ F[[\g^*]] $ et soit  $ \, \e \, $  l'id\'eal maximal
(unique) de  $ \, F[[\g^* \oplus \g^*]] = F[[\g^*]] \otimes F[[\g^*]]
\, $  (produit tensoriel topologique, selon [Di], Ch.~1).  Puisque
$ {\frak R} $  est un automorphisme d'alg\`ebre,  $ \, {\frak R}(\e)
= \e \, $  et  $ {\frak R} $  induit un automorphisme
$ \overline{\frak R} $  de l'espace vectoriel  $ \e \big/ \e^2 $.
Or  $ \, \e \big/ \e^2 $  s'identifie \`a l'alg\`ebre de Lie
$ \g \oplus \g \, $  et comme $ {\frak R} $  est un automorphisme
d'alg\`ebre de Poisson, l'application  $ \overline{\frak R} $  est
un automorphisme de l'alg\`ebre de Lie  $ \, \g \oplus \g $;
l'automorphisme  $ \overline{\frak R} $  h\'erite aussi des
autres propri\'et\'es du tressage  $ {\frak R} $.  En particulier,
$ {\frak R} $  et  $ \overline{\frak R} $  sont solutions de la QYBE
et d\'efinissent donc une action du groupe des tresses  $ {\Cal B}_n $
sur  $ {F[[\g^* \oplus \g^*]]}^{\otimes n} $  et sur
$ {(\g \oplus \g)}^{\otimes n} $.
                                              \par
   De tels automorphismes de  $ \, \g \oplus \g \, $  ont \'et\'es
introduits dans [WX], \S 9; leur construction est li\'ee \`a la
"$ R $-matrice  globale", qui donne aussi une interpr\'etation
g\'eom\'etrique de la  $ r $-matrice  classique.  Il conviendrait
de comparer nos r\'esultats et ceux de [WX] et d'\'etudier
parall\`element les propri\'et\'es de fonctorialit\'e de
notre construction.

\vskip1,1truecm

\centerline {\bf \S \; 3. \  D\'emonstration du th\'eor\`eme 2.1 }

\vskip10pt

   Dans cette section  $ (H,R) $  sera une alg\`ebre de Hopf
quasitriangulaire comme dans l'\'enonc\'e du th\'eor\`eme 2.1.  Nous
voulons \'etudier l'action adjointe de  $ R $  sur
l'alg\`ebre  $ H \otimes H $. Cette
derni\`ere poss\`ede une structure naturelle d'alg\`ebre de Hopf,
son coproduit $ \tilde\Delta $ \'etant d\'efini par  $ \;
\tilde\Delta = \sigma_{2{}3} \circ (\Delta \otimes
{\Id}_{\scriptscriptstyle H} \otimes {\Id}_{\scriptscriptstyle H})
\circ ({\Id}_{\scriptscriptstyle H} \otimes \Delta) $,
o\`u  $ \, \sigma_{2{}3} \, $  d\'esigne la volte dans les positions
$ 2 $  et  $ 3 $.  Nous noterons aussi  $ \, I = 1 \otimes 1 \, $
l'unit\'e dans  $ H \otimes H $.  Selon notre d\'efinition du produit
tensoriel dans  $ {\Cal A} $,  on a  $ \, {\big( H \otimes H \big)}' =
H' \otimes H' \, $.  Notre but est de montrer que, m\^eme si $ R $
n'appartient pas \`a  $ {\big( H \otimes H \big)}' $,  son action
adjointe  $ \, a \mapsto R \cdot a \cdot R^{-1} \, $  laisse
stable  $ \, {\big( H \otimes H \big)}' = H' \otimes H' \, $.
                                                  \par
   Posons tout d'abord, pour  $ \, \Sigma = \{i_1, \dots, i_k\} \subseteq
\{1, \dots,n\} \, $,  toujours avec  $ \, i_1 < \cdots < i_k \, $:
  $$  R_\Sigma = R_{2i_1-1,2i_k} R_{2i_1-1,2i_{k-1}} \cdots
R_{2i_1-1,2i_1} R_{2i_2-1,2i_k} \cdots R_{2i_{k-1}-1,2i_1} R_{2i_k-1,2i_k}
\cdots R_{2i_k-1,2i_1}  $$
(produit de  $ k^2 $  termes) o\`u  $ \, R_{r,s} =
j_{\scriptscriptstyle \{r,s\}}(R) \, $,  en d\'efinissant  $ \,
j_{\scriptscriptstyle \{r,s\}} \colon \, H \otimes H \longrightarrow
H^{\otimes 2n} \, $  comme pr\'ec\'edemment.  Nous noterons toujours
$ |\Sigma| $  pour le cardinal de  $ \Sigma $  (ici  $ \, |\Sigma| =
k \, $).

\vskip7pt

\proclaim{Lemme 3.1}
   Dans  $ {\big( H \otimes H \big)}^{\otimes n} $,  pour tout  $ \,
\Sigma \subseteq \{1, \dots, n\} $,  on a  $ \, {\tilde\Delta}_\Sigma(R) =
R_\Sigma \, $.
\endproclaim

\demo{Preuve}  Sans nuire \`a la g\'en\'eralit\'e du probl\`eme, nous pouvons nous contenter
de prouver le
r\'esultat pour  $ \, \Sigma \! = \! \{1, \dots, n\} $,  {\it i.e.},
  $$  {\tilde\Delta}_{\{1, \dots, n\}}(R) = R_{\{1, \dots, n\}} = R_{1,2n}
\cdot R_{1,2n-2} \cdots R_{1,2} \cdot R_{3,2n} \cdots R_{2n-3,2} \cdot
R_{2n-1,2n} \cdots R_{2n-1,2} \; .  $$
   \indent   Le r\'esultat est \'evident au rang  $ \, n = 1 \, $.
Supposons-le acquis au rang  $ \, n\geq 1 \, $,  et montrons-le au rang  $ \, n
+ 1 \, $:  par  d\'efinition de  $ \tilde\Delta $  et par les
propri\'et\'es de la  $ R $-matrice  on a
  $$  \eqalign{
   {\tilde\Delta}_{\{1, \dots, n+1 \}}  (\!R)  &  = \left( {\tilde\Delta}
\otimes {{\Id}^{\, \otimes {n-1}}_{\scriptscriptstyle H \otimes H}}
\right) \! \big( {\tilde\Delta}_{\{1, \dots, n\}} (R) \big)\cr
 & = \left(
{\tilde\Delta} \otimes {{\Id}^{\,
\otimes {n-1}}_{\scriptscriptstyle H \otimes H}} \right) \! \big( R_{\{1, \dots, n\}} \big)  \cr
   &  = \sigma_{2{}3} (\Delta \otimes {\Id}_{\scriptscriptstyle H}^{\,
\otimes  2n} )\!\!\left( {\Id}_{\scriptscriptstyle H} \otimes \Delta
\otimes  {\Id}_{\scriptscriptstyle H}^{ \otimes 2(n-1)}\! \right) \!\!
(R_{1,2n} \cdots  \hskip-0,16pt  R_{1,2} \cdots  \hskip-0,16pt  R_{3,2}
\cdots  \hskip-0,16pt   R_{2n-1,2})  \cr
   &  = \sigma_{2{}3} \! \left( \Delta \otimes {\Id}_{\scriptscriptstyle
H}^{\, \otimes  2n} \right) (R_{1,2n+1} \cdots R_{1,3} R_{1,2} \cdots
R_{4,3} R_{4,2} \cdots  R_{2n,3} R_{2n,2})  \cr
   &  = \sigma_{2{}3} (R_{1,2n+2} R_{2,2n+2} \cdots R_{1,4} R_{2,4}
R_{1,3}  R_{2,3} \cdots R_{5,4} R_{5,3} \cdots R_{2n+1,4} R_{2n+1,3})
\cr
   &  = R_{1,2n+2} R_{3,2n+2} \cdots R_{1,4} R_{3,4} \cdot R_{1,2} R_{3,2}
\cdots R_{5,4} \cdot R_{5,2} \cdots R_{2n+1,4} R_{2n+1,2}  \cr
   &  = R_{1,2n+2} \cdots R_{1,4} R_{1,2} R_{3,2n+2} \cdots R_{3,4}
R_{3,2} \cdots R_{5,4} R_{5,2} \cdots R_{2n+1,4} R_{2n+1,2}  \cr
   &  = R_{\{1, \dots, n+1\}}.\quad
\square  \cr }  $$
\enddemo

\vskip7pt

Remarque: dor\'enavant pour tous
$ \, a $,  $ b \in \N \, $,  nous utiliserons la notation
$ \, C^a_b \, $  pour d\'esigner l'entier  $ \, {b \choose a}
= {\, b! \, \over \, a! (b-a)! \,}$.

\vskip7pt

\proclaim{Lemme 3.2}  Pour tout  $ \, a \in {\big( H \otimes
H \big)}'$  et pour tout ensemble  $ \, \Sigma \, $  tel que
$ \, |\Sigma| > i \, $,  on a
  $$  {\tilde\Delta}_\Sigma(a) = \sum_{\Sigma' \subseteq \Sigma, \;\,
|\Sigma'| \leq i}  \!\! {(-1)}^{i-|\Sigma'|} \, C_{|\Sigma| - 1 -
|\Sigma'|}^{i-|\Sigma'|} \, {\tilde\Delta}_{\Sigma'}(a) + O \big(
h^{i+1} \big) \, .  $$
\endproclaim

\demo{Preuve}  Il suffit de prouver l'\'enonc\'e pour  $ \, \Sigma = \{1,
\dots, n\} , $  avec  $ \, n > i \, $.  Gr\^ace \`a (1.2), on a
  $$  \eqalign{
   {\tilde\Delta}_{\{1, \dots, n\}}(a)  &  = \sum_{\bar\Sigma \subseteq
\{1, \dots, n\}} \delta_{\bar\Sigma}(a) \cr
& = \sum_{\bar\Sigma \subseteq \{1,
\dots,  n\},~
|\bar\Sigma| \leq i } \! \delta_{\bar\Sigma}(a) + O \big(
h^{i+1} \big)  \cr
   &  = \sum_{\Sb \bar\Sigma \subseteq \{1, \dots, n\}\\
   |\bar\Sigma| \leq i \endSb } \;
\sum_{\Sigma' \subseteq \bar\Sigma} {(-1)}^{|\bar\Sigma|-|\Sigma'|} \,
{\tilde\Delta}_{\Sigma'}(a) + O \big( h^{i+1} \big)  \cr
   &  = \sum_{\Sb \Sigma' \subseteq \{1, \dots, n\}\\
   |\Sigma'| \leq i\endSb } \!\!\!
{\tilde\Delta}_{\Sigma'}(a) \sum_{\Sigma' \subseteq
\bar\Sigma,~|\bar\Sigma|
\leq i} {(-1)}^{|\bar\Sigma|-|\Sigma'|} + O \big( h^{i+1} \big)  \cr
   &  = \sum_{\Sigma' \subseteq \{1, \dots, n\},~
   |\Sigma'| \leq i} \!\!\!
{\tilde\Delta}_{\Sigma'}(a) \, {(-1)}^{i-|\Sigma'|} \,
C^{i-|\Sigma'|}_{n-1-|\Sigma'|} + O \big( h^{i+1} \big).  \;\;
 \quad  \square  \cr }  $$
\enddemo

\vskip7pt

   Avant de nous attaquer au r\'esultat principal, il nous faut encore un
petit rappel technique sur les coefficients du bin\^ome qui peut se
prouver facilement en utilisant le d\'eveloppement en s\'erie formelle de
$ \, {(1-X)}^{-(r+1)} \, $,  \`a savoir  $ \, {(1-X)}^{-(r+1)} =
\sum\limits_{k=0}^{\infty} C_{k+r}^r  X^k \, $.

\vskip7pt

\proclaim{Lemme 3.3}  Soient  $ \, r $,  $ s $,  $ t \in \N \, $  tels que
$ \, r < t $.  On a alors les relations suivantes (o\`u l'on pose  $ \,
C_u^v = 0 \, $  si  $ \, v > u \, $):
                            \hfill\break
   \centerline{ $ \displaystyle{ (a) \quad  \sum_{d=0}^{t} {(-1)}^{d}
\, C_{d-1}^{r} \, C_{t}^{d}  = -{(-1)}^{r} \; ,  \qquad   (b) \quad
\sum_{d=0}^{t} {(-1)}^{d} \, C_{d+s}^{r} \, C_{t}^{d}  = 0 .
 } $ }
\endproclaim

   Voici enfin le r\'esultat principal de cette section:

\proclaim{Proposition 3.4} On a $ \; R \, a \, R^{-1} \in {\big( H \otimes H
\big)}' \, $ pour tout  $ \, a \in {\big( H \otimes H
\big)}' \, $.
\endproclaim

\demo{Preuve}  Comme nous devons montrer que  $ \, R \, a \, R^{-1} \, $
appartient \`a  $ {\big( H \otimes H \big)}'$,  nous devons
consid\'erer les  termes  $ \delta_n \! \left( R \, a \, R^{-1} \right) $,
$ \, n \in \N \, $.  Pour cela r\'e\'ecrivons  $ \, \delta_{\{1, \dots,
n\}} \! \left( R \, a \, R^{-1} \right) $  en utilisant le lemme 3.1 et
le fait que  $  \tilde\Delta $,  et plus g\'en\'eralement
$ {\tilde\Delta}_{\{i_1, \dots, i_k\}} $  (pour  $ k \leq n $),  est un
morphisme d'alg\`ebre:  $ \; \delta_{\{1, \dots, n\}} \left( R \, a \,
R^{-1} \right) = \sum\limits_{\Sigma \subseteq \{1, \dots, n\}}
{(-1)}^{n-|\Sigma|} R_\Sigma \, {\tilde\Delta}_\Sigma(a) \,
R^{-1}_\Sigma \; $.
                                                \par
   Nous allons d\'emontrer par r\'ecurrence sur  $ i $  que
  $$  \delta_{\{1, \dots, n\}} \left( R \, a \, R^{-1} \right) = O \big(
h^{i+1} \big)  \qquad  \text{pour tout} \quad  0 \leq i \leq n-1 \; .
\eqno (\star)  $$
Nous verrons ainsi que tous les termes du d\'eveloppement limit\'e
\`a  l'ordre  $ \, n-1 \, $  sont nuls, donc que $ \, \delta_n \left( R \, a
\, R^{-1} \right) = O(h^n) \, $,  d'o\`u notre \'enonc\'e.
                                                \par
   Pour  $i = 0$ et  pour chaque  $ \Sigma \, $, on a
$ \, {\tilde\Delta}_\Sigma(a) = \epsilon(a) I^{\otimes n} + O(h) \, $,
$ \, R_\Sigma = I^{\otimes n} + O(h) \, $,  et aussi
$ \, R_\Sigma^{-1} = I^{\otimes n} + O(h) \, $,  d'o\`u
$$ \;
\delta_{\{1, \dots, n\}} \left( R \, a \, R^{-1} \right) =
\sum\limits_{k=1}^{n} C_n^k {(-1)}^{n-k} \epsilon(a) \, I^{\otimes n}
+ O(h) = O(h) \;,$$
donc le r\'esultat  $ (\star) $  est vrai pour
$ \, i = 0 \, $.
                                                \par
   Supposons le r\'esultat  $ (\star) $  acquis pour tout  $ \, i' < i
\, $.  \'Ecrivons les d\'eveloppements  $ h $-adiques  de
$ R_\Sigma $  et  $ R_\Sigma^{-1} $  sous la forme  $ \; R_\Sigma
= \sum_{\ell=0}^\infty R_\Sigma^{\,(\ell)} \, h^\ell \; $  et  $ \;
R_\Sigma^{-1} = \sum_{m=0}^\infty R_\Sigma^{\,(-m)} \, h^m \, $.
Le lemme 3.2 fournit une approximation de
$ {\tilde\Delta}_\Sigma(a) $  \`a l'ordre  $ j $:
  $$  {\tilde\Delta}_\Sigma(a) \, = \sum_{\Sigma' \subseteq
\Sigma,~|\Sigma'|\leq j} {(-1)}^{j-|\Sigma'|} \, C_{|\Sigma| - 1 -
|\Sigma'|}^{j  - |\Sigma'|} \, {\tilde\Delta}_{\Sigma'}(a) + O \big(
h^{j+1} \big) \; .  $$
   Nous obtenons l'approximation suivante de $ \, \delta_{\{1, \dots, n\}}
\left( R \, a \, R^{-1} \right) $:
  $$  \displaylines{
   {} \;   \delta_{\{1, \dots, n\}} \left( R \, a \, R^{-1} \right) =
\sum_{\Sigma \subseteq \{1, \dots, n\}} \sum_{\ell + m \leq i} {(-1)}^{n -
|\Sigma|} \, R_\Sigma^{\,(\ell)} \, {\tilde\Delta}_\Sigma (a) \,
R_\Sigma^{\,(-m)} \, h^{\ell + m} + O \big( h^{i+1} \big) =   \hfill {\ }
\cr
   {} \,   = \sum_{j=0}^i \,\, \sum_{\ell + m = i - j} \Bigg(
\sum_{\Sb  \Sigma \subseteq \{1, \dots, n\}  \\  |\Sigma| > j  \endSb}
\sum_{\Sb  \Sigma' \subseteq \Sigma  \\  |\Sigma'| \leq j  \endSb}  \!
{(-1)}^{n-|\Sigma|} \, {(-1)}^{j-|\Sigma'|} \,
C_{|\Sigma|-1-|\Sigma'|}^{j-|\Sigma'|} \, R^{\,(\ell)}_\Sigma \,
{\tilde\Delta}_{\Sigma'}(a) \, R^{\,(-m)}_\Sigma +   \hfill {\ }  \cr
   {\ } \hfill   + \sum_{\Sb  \Sigma \subseteq \{1, \dots, n\}  \\
|\Sigma| \leq j  \endSb}  {(-1)}^{n-|\Sigma|} \, R^{\,(\ell)}_\Sigma \,
{\tilde\Delta}_\Sigma(a) \, R^{\,(-m)}_\Sigma \Bigg) \, h^{\ell + m}
+ O \big( h^{i+1} \big) =  \cr
   {} \; \hfill   = \sum_{j=0}^i \, \sum_{\ell +m+j=i} \, \sum_{\Sb
\Sigma' \subseteq \{1, \dots, n\}  \\  |\Sigma'| \leq  j  \endSb}  \!\!
\Bigg( \sum_{\Sb  \Sigma \subseteq \{1, \dots, n\}  \\  \Sigma' \subseteq
\Sigma, {\ }  |\Sigma| > j  \endSb}  \hskip-14pt  {(-1)}^{n-|\Sigma|} \,
{(-1)}^{j-|\Sigma'|} \, C_{|\Sigma|-1-|\Sigma'|}^{j-|\Sigma'|} \,
R^{\,(\ell)}_\Sigma {\tilde\Delta}_{\Sigma'}(a) \, R^{\,(-m)}_\Sigma +
{\ }  \cr
   {\ } \hfill   + {(-1)}^{n-|\Sigma'|} \, R^{\,(\ell)}_{\Sigma'} \,
{\tilde\Delta}_{\Sigma'}(a) \, R^{\,(-m)}_{\Sigma'} \Bigg) \, h^{\ell + m}
+ O \big( h^{i+1} \big) \; .  \cr }  $$
   \indent   Notons (E) la derni\`ere expression entre
parenth\`eses; nous montrerons que cette expression est nulle, d'o\`u
$ \, \delta_n \left( R \, a \, R^{-1} \right) = O \big( h^{i+1} \big)
\, $.

\medskip

   Regardons d'abord les termes correspondant \`a  $ \, \ell + m = 0
\, $,  c'est-\`a-dire  $ \, j = i \, $.  Nous retrouvons $ \delta_{\{1,
\dots, n\}}(a) $,  qui est dans  $ O \big( h^{i+1} \big) $  par
hypoth\`ese.  Dans la suite du calcul nous supposerons d\'esormais
$ \, \ell + m > 0 \, $.

   Fixons maintenant un entier strictement positif $S$ et
regardons  comment les termes  $ R_\Sigma^{\,(\ell)} $  et
$ R_\Sigma^{\,(-m)} $  agissent sur  $ {{\big( H \otimes H \big)}'}^{\,
\otimes n} $  (respectivement \`a gauche et \`a droite) pour  $ \, \ell +
m =S\, $.  En faisant
le d\'eveloppement limit\'e de chaque  $ R_{i,j} $  qui appara\^\i{}t
dans  $ R_\Sigma \, $,  on voit que  $ R_\Sigma^{\,(\ell)} $  et
$ R_\Sigma^{\,(-m)} $  sont sommes de produits d'au plus  $ \ell $
et  $ m $  termes respectivement, chacun agissant sur au plus deux
facteurs tensoriels de  $ {{\big( H \otimes H \big)}'}^{\, \otimes n} $.
Nous allons r\'e\'ecrire  $ \, \sum\limits_{\ell + m = S}
R^{\,(\ell)}_\Sigma \, {\tilde\Delta}_{\Sigma'}(a) \,
R^{\,(-m)}_\Sigma \, $  en regroupant les termes de la somme
qui agissent sur les m\^emes facteurs de  $ {{\big( H \otimes H
\big)}'}^{\, \otimes n} $,  facteurs dont nous identifierons les
positions par  $ \, \Sigma'' $.
                                             \par
   Si  $ i $  appartient \`a  $ \Sigma'' $,  dans l'identification
$ \, {(H \otimes H)}^{\otimes n} = H^{\otimes 2n} \, $  que
nous avons choisie pour d\'efinir  $ R_\Sigma \, $,  l'indice  $ i $
correspond \`a la paire  $ (2i-1,2i) \, $;  mais alors  $ R_\Sigma $  et
$ R_\Sigma^{\,-1} $,  et  donc aussi chaque  $ R_\Sigma^{\,(\ell)} $  et
chaque  $ R_\Sigma^{\,(-m)} \, $,  n'agissent de mani\`ere non triviale sur
le  $ i $-\`eme  facteur de  $ {\tilde\Delta}_{\Sigma'}(a) $  que si, dans
l'\'ecriture explicite de  $ R_\Sigma $, un terme non trivial apparait aux
places  $ \, 2i-1 \, $  ou  $ 2i $,  donc seulement si  $ \, i \in \Sigma
\, $:  ainsi  $ \, \Sigma'' \subseteq \Sigma \, $.  Nous posons alors
  $$  \sum_{\ell + m = S} R^{\,(\ell)}_\Sigma \,
{\tilde\Delta}_{\Sigma'}(a) \, R^{\,(-m)}_\Sigma =
\sum_{\Sigma'' \subseteq \Sigma} A^{(S)}_{\Sigma', \Sigma,
\Sigma''}(a) \, .  $$
   \indent   Maintenant consid\'erons  $ \, \bar\Sigma \supseteq \Sigma
\, $.  D'apr\`es la d\'efinition on a  $ \, R_{\bar\Sigma} = R_\Sigma +
{\Cal A} \, $,  o\`u  $ {\Cal A} $  est une somme de termes qui
contiennent des facteurs  $ \, R_{2i-1,2j}^{\,(s)} \, $  avec  $ \, \{i,j\}
\not\subseteq \Sigma \, $:  pour d\'emontrer ceci, il suffit de d\'evelopper chaque
facteur  $ R_{a,b} $  dans  $ R_{\bar\Sigma} $  comme  $ \, R_{a,b} =
1^{\otimes 2n} + O(h) \, $.  De  m\^eme, on a aussi  $ \,
R_{\bar\Sigma}^{\,(\ell)} = R_\Sigma^{\,(\ell)} + {\Cal A}' \, $,  et
pareillement  $ \, R_{\bar\Sigma}^{\,(-m)} = R_\Sigma^{\,(-m)} + {\Cal A}''
\, $. Cela implique que  $ \,  A^{(S)}_{\Sigma'', \bar\Sigma, \Sigma'}(a) =
A^{(S)}_{\Sigma'', \Sigma, \Sigma'}(a) \, $,  et donc que les
$ A^{(S)}_{\Sigma'', \Sigma, \Sigma'}(a) $  ne d\'ependent pas de
$ \Sigma \, $;  on \'ecrit alors
  $$  \sum_{\ell + m = S} R^{\,(\ell)}_\Sigma \,
{\tilde\Delta}_{\Sigma'}(a) \,  R^{\,(-m)}_\Sigma = \sum_{\Sigma''
\subseteq \Sigma}  A^{(S)}_{\Sigma',\Sigma''}(a) \; .  $$
   \indent   Nous allons ensuite r\'e\'ecrire  $ (E) $  \`a l'aide des
$ A^{(S)}_{\Sigma', \Sigma''}(a) $.  Par commodit\'e dans la suite des
calculs, nous noterons  $ \, \delta_{\Sigma'' \subseteq \Sigma'} \, $  la
fonction  qui vaut  $ 1 $  si  $ \, \Sigma'' \subseteq \Sigma' \, $
et  $ 0 $  sinon.  Nous obtenons alors une nouvelle expression pour
$ \, \delta_{\{1, \dots,  n\}} \left( R \, a \, R^{-1} \right) \, $,
\`a savoir
  $$  \displaylines{
   {} \;   \delta_{\{1, \dots, n\}} \! \left( R \, a \, R^{-1} \right)
= \sum_{j=0}^{i-1} \, \sum_{\Sb  \Sigma' \subseteq \{1, \dots, n\}  \\
   |\Sigma'| \leq  j  \endSb} \Bigg( \sum_{\Sb  \Sigma \subseteq \{1,
\dots, n\}  \\  \Sigma' \subseteq \Sigma, {\ } |\Sigma| > j  \endSb}
\!\!\!\! {(-1)}^{n-|\Sigma|} \, {(-1)}^{j-|\Sigma'|} \,
C_{|\Sigma|-1-|\Sigma'|}^{j-|\Sigma'|} \times  \cr
   {\ } \hfill   \times \sum_{\Sigma'' \subseteq \Sigma}
\! A^{(i-j)}_{\Sigma', \Sigma''}(a) + {(-1)}^{n-|\Sigma'|}
\sum_{\Sigma'' \subseteq \Sigma'} A^{(i-j)}_{\Sigma', \Sigma''}(a)
\Bigg) \, h^{i-j} + O \big( h^{i+1} \big) {\ }  \cr
   {\;} = \, \sum_{j=0}^{i-1} \, \sum_{\Sb  \Sigma' \subseteq
\{1, \dots, n\}  \\  |\Sigma'| \leq  j  \endSb}  h^{i-j}  \sum_{\Sigma''
\subseteq \{1, \dots, n\}} \! A^{(i-j)}_{\Sigma', \Sigma''}(a) \times
\hfill {\ }  \cr
   {\ } \hfill   \times \Bigg( \sum_{\Sb  \Sigma \subseteq \{1, \dots,
n\}  \\
\Sigma' \subseteq \Sigma, \; \Sigma'' \subseteq \Sigma, \; |\Sigma| > j
\endSb}  \!\!\!\!\! {(-1)}^{n-|\Sigma|} \, {(-1)}^{j-|\Sigma'|} \,
C_{|\Sigma| -  1 - |\Sigma'|}^{j - |\Sigma'|} + {(-1)}^{n-|\Sigma'|} \,
\delta_{\Sigma'' \subseteq \Sigma'} \Bigg)  + O \big( h^{i+1} \big) \; .
{\ }  \cr }  $$
   \indent   Notons  $ {\big( E' \big)}_{\Sigma',\Sigma''} $
la nouvelle expression entre parenth\`ese; autrement dit, pour
$ \Sigma' $  et  $ \Sigma'' $  fix\'ees, avec  $ \, \big| \Sigma' \big|
\leq j \, $,  on pose
  $$  {\big( E' \big)}_{\Sigma',\Sigma''} = \sum_{\Sb  \Sigma \subseteq
\{1, \dots, n\}  \\  \Sigma' \subseteq \Sigma, \; \Sigma'' \subseteq
\Sigma, \; |\Sigma| > j  \endSb} \!\!\!\!\! {(-1)}^{n-|\Sigma|} \,
{(-1)}^{j-|\Sigma'|} \, C_{|\Sigma| - 1 -  |\Sigma'|}^{j - |\Sigma'|} +
{(-1)}^{n-|\Sigma'|} \, \delta_{\Sigma'' \subseteq \Sigma'}. $$
Remarquons que cette expression est purement
combinatoire. Nous allons d\'emontrer qu'elle est nulle lorsque les parties
$ \Sigma' $  et  $ \Sigma'' $  sont telles que  $ \, \big| \Sigma' \cup
\Sigma'' \big| \leq i - j + \big| \Sigma' \big| \, $  et  $\,\big| \Sigma'
\big| \leq j \,$.  En vertu du lemme suivant, ceci suffira pour
prouver la proposition. Remarquons d\'ej\`a que
$ \, j < i \, $  et
$ \, i \leq n-1 \, $,  donc  $ \, j \leq n-2$.

\vskip7pt

\proclaim{Lemme 3.5}
\par
Pour tout  $ \, S > 0 \, $,  dans
l'expression  $ \; \sum\limits_{\ell + m = S} R^{\,(\ell)}_\Sigma \,
{\tilde\Delta}_{\Sigma'}(a) \, R^{\,(-m)}_\Sigma = \sum\limits_{\Sigma''
\subseteq \Sigma} A^{(S)}_{\Sigma',\Sigma''}(a) \, $  \; on a
$ \; A^{(S)}_{\Sigma', \Sigma''}(a) = 0 $  \ pour tout  $ \Sigma' $,
$ \Sigma'' $  tels que  $ \, \big| \Sigma' \cup \Sigma'' \big| > S +
\big| \Sigma' \big| \, $.
\endproclaim

\demo{Preuve}  Pour d\'emontrer ce r\'esultat
nous \'etudions l'action adjointe de  $ R_\Sigma $  sur
$ {\big( H \otimes H \big)}^{\, \otimes n} $.
                                          \par
   Premi\`erement, sur  $ k \cdot I^{\otimes n} $  l'action de ces
\'el\'ements  donne un terme nul car pour  $ \, S > 0 \, $,
on retrouve le terme \`a l'ordre  $ S
$  du  d\'eveloppement  $ h $-adique  de  $ \, R_\Sigma \cdot
R_\Sigma^{-1} = 1$.
                                          \par
   Deuxi\`emement, consid\'erons $ \, \Sigma \subseteq \{1, \dots, n\}
\, $,   et \'etudions l'action sur
$$ \, {\big( H \otimes H)}_{\Sigma'} =
j_{\scriptscriptstyle \Sigma'} \left( {\big( H \otimes H \big)}^{\otimes
|\Sigma'|} \right) \, \subseteq {\big( H \otimes H \big)}^{\otimes n}. $$
$ R_\Sigma $  est un produit de  $ {|\Sigma|}^2 $
termes du  type  $R_{a,b}$,  avec  $a, b \in \big\{\, 2i-1, 2j
\,\big\vert\, i,  j \in \Sigma \,\big\}$;  analysons ce qui se
passe lorsqu'on effectue le  produit  $ \, P = R_\Sigma \cdot x \cdot
R_\Sigma^{\,-1}$  si  $ \, x \in  {\big( H \otimes H)}_{\Sigma'} \, $.
                                         \par
   Consid\'erons le facteur  $ R_{a,b} $  qui appara\^\i{}t le plus \`a
droite: si  $ \, a, b \not\in \big\{\, 2j-1, 2j \,\big\vert\, j \in
\Sigma' \,\big\} \, $,  alors en calculant  $ P $  on trouve  $ \, P =
R_\Sigma \, x \,  R_\Sigma^{\,-1} = R_\star \, R_{a,b} \, x \,
R_{a,b}^{\,-1} \, R_\star^{\,-1} =  R_\star \, x \, R_\star^{\,-1} \, $
(o\`u  $ \, R_\star = R_\Sigma \,  R_{a,b}^{\,-1} \, $).  De m\^eme, en
avan\c{c}ant de droite \`a gauche le long  de  $ R_\Sigma $  on peut
\'ecarter tous les facteurs  $ R_{c,d} $  de ce type, \`a savoir
les facteurs tels que
$ \, c, d \not\in \big\{\, 2j-1, 2j \,\big\vert\, j \in \Sigma' \,\big\}
\, $.  Ainsi le premier facteur dont l'action adjointe est non  triviale
sera n\'ecessairement du type  $ R_{\bar{a},\bar{b}} $  avec l'un des
deux indices appartenant \`a  $ \, \big\{\, 2j-1, 2j \,\big\vert\, j
\in \Sigma' \,\big\} \, $,  soit par exemple  $ \bar{a} $.  Notons que le
nouvel indice  $ \, \bar{a} \; (\, \in \{1, 2, \dots, 2n-1, 2n\} \,)
\, $,  qui  agit sur un facteur tensoriel dans  $ H^{\otimes 2n} $,
correspond \`a un nouvel indice  $ \, j_{\bar{a}} \; (\, \in \{1, \dots,
n\} \,) \, $,  agissant  sur un facteur tensoriel de  $ {\big( H \otimes
H \big)}^{\otimes n} \, $.  Ainsi pour les facteurs successifs   ---
c'est-\`a-dire \`a gauche de  $ R_{\bar{a},\bar{b}} $  ---   il faut r\'ep\'eter la
m\^eme analyse, mais avec l'ensemble  $ \, \big\{\, 2j-1, 2j \,\big\vert\,
j \in \Sigma' \cup \{ j_{\bar{a}} \} \,\big\} \, $  \`a la place de  $ \,
\big\{\, 2j-1, 2j \,\big\vert\, j \in \Sigma' \,\big\} \, $;  donc, comme
$ R_{\bar{a},\bar{b}} $ ne  peut agir de mani\`ere non triviale que sur au
plus  $ \big| \Sigma' \big| $  facteurs de  $ {\big( H \otimes H
\big)}^{\otimes n} $, le facteur \`a sa gauche ne
peut agir de mani\`ere non triviale que sur au plus  $ \, \big| \Sigma'
\big| + 1 \, $  facteurs.  La conclusion est que  l'action adjointe
de  $ R_\Sigma $  est non triviale sur au plus  $ \, \big| \Sigma'
\big| + \big| \Sigma \big| \, $  facteurs de
$ {\big( H \otimes H \big)}^{\otimes n} $.
                                          \par
   Maintenant, consid\'erons les diff\'erents termes
$ R_\Sigma^{\,(\ell)} $   et  $ R_\Sigma^{\,(-m)} \,$,
avec  $ \, \ell + m = S \, $,  et \'etudions les  produits
$ \, R_\Sigma^{\,(\ell)} \cdot x \cdot R_\Sigma^{\,(-m)} \, $,
avec  $ \, x \in {\big( H \otimes H)}_{\Sigma'} \, $.  On sait d\'ej\`a que
$ R_\Sigma^{\,(\ell)} $  et  $ R_\Sigma^{\,(-m)} \, $  sont sommes de
produits,  not\'es  $ P_+ $  et  $ P_- \, $,  d'au plus  $ \ell $  et
$ m $  termes  respectivement, du type  $ \, R_{i,j}^{\,(\pm k)} \, $;
les termes  $ A_{\Sigma',\Sigma''}^{(S)}(a) $  ne sont alors que des sommes
de termes du type  $ \, P_+ \, {\tilde\Delta}_{\Sigma'}(a) \, P_- \, $,
o\`u de plus les ``indices'' intervenant dans  $ P_+ $  et  $ P_- \, $
sont dans  $ \Sigma'' $.  Or,  comme chaque  $ P_+ $  et  chaque  $ P_- $
est un produit d'au plus  $ \ell $   et  $ m $  facteurs  $ \,
R_{i,j}^{\,(\pm k)} \, $,  on peut raffiner l'argument  pr\'ec\'edent.
Consid\'erons seulement le terme \`a l'ordre  $ S $  du d\'eveloppement
$ h $-adique  de  $ \, P = R_\Sigma \, x \,  R_\Sigma^{\,-1} = R_\star
\, R_{a,b} \, x \, R_{a,b}^{\,-1} \, R_\star^{\,-1} =  R_\star \, x \,
R_\star^{\,-1} \, $:  lorsqu'il y a des facteurs du type
$ R_{a,b}^{\,(k)} $  ou  $ R_{a,b}^{\,(t)} \, $,  pour  $ a $ et  $ b $
fix\'es  n'appartenant pas \`a  $ \big\{\, 2j-1, 2j \,\big\vert\,
j \in \Sigma' \,\big\} \, $, qui apparaissent dans
$ R_\Sigma^{\,(\ell)} $   ou  $ R_\Sigma^{\,(-m)} \, $ pour certains
$ \ell $  ou  $ m $,  la contribution totale de tous ces termes dans la
somme  $ \, \sum\limits_{\ell + m = S} R_\Sigma^{\,(\ell)} \, x \,
R_\Sigma^{\,(-m)} \, $  est nulle.  De plus, comme on ne consid\`ere
maintenant que  $ S $  facteurs, on conclut que  $ \; A_{\Sigma',
\Sigma''}^{(S)}(a) = 0 \; $  si  $ \; \big| \Sigma' \cup \Sigma''
\big| > S + \big| \Sigma' \big| \, $.   $ \square $
\enddemo

\vskip5pt

   Revenons \`a la d\'emonstration de la proposition 3.4 et
calculons $ {\big( E' \big)}_{\Sigma',\Sigma''} $.
Gr\^ace \`a la remarque pr\'ec\'edente, nous pouvons nous limiter aux
paires  $ \big( \Sigma', \Sigma'' \big) $  telles que
$$ \; \big| \Sigma'
\cup \Sigma'' \big| \leq i - j + \big| \Sigma' \big| \leq i - j + j =
i \leq n-1 \;.$$
On pourra toujours trouver au moins deux parties $
\Sigma $ de $ \{1, \dots, n\} \, $  telles que  $ \, |\Sigma| > j \, $
et  $ \, \Sigma' \cup \Sigma'' \subseteq \Sigma \, $,  ce qui nous assure
qu'il y aura toujours au  moins deux termes dans le comptage qui va suivre
(condition qui assurera la  nullit\'e de l'expression  $ {\big( E'
\big)}_{\Sigma',\Sigma''} \, $).  Nous  allons distinguer trois cas:

\vskip3pt

 (I) \quad  Si  $ \, \Sigma'' \subseteq \Sigma' \, $,  alors l'expression
$ {\big( E' \big)}_{\Sigma',\Sigma''} $  devient
  $$  {\big( E':1 \big)}_{\Sigma',\Sigma''}
 \, = \sum_{\Sb  \Sigma \subseteq \{1, \dots, n\}  \\
\Sigma' \subseteq \Sigma, \; |\Sigma| > j  \endSb}  \!\!
{(-1)}^{n-|\Sigma|} \,  {(-1)}^{j-|\Sigma'|} \,
C_{|\Sigma|-1-|\Sigma'|}^{j-|\Sigma'|} +  {(-1)}^{n-|\Sigma'|} \; .  $$
En regroupant les  $ \Sigma $  qui ont le m\^eme cardinal  $ d $,  un
simple  comptage nous donne
  $$  { \big( E':1 \big)}_{\Sigma',\Sigma''} = \sum_{d=j+1}^n {(-1)}^{n-d}
\,  {(-1)}^{j-|\Sigma'|}
\, C_{d-1-|\Sigma'|}^{j-|\Sigma'|} \, C_{n-|\Sigma'|}^{d-|\Sigma'|} +
{(-1)}^{n-|\Sigma'|} \; .  $$
   \indent   Or, cette derni\`ere expression est nulle d'apr\`es le
lemme 3.3, car elle correspond \`a une somme du type
$$ \;
\sum\limits_{k=r+1}^t {(-1)}^{t+r-k} \, C_{k-1}^r \, C_t^k + {(-1)}^t =
\sum\limits_{k=0}^t {(-1)}^{t+r-k} \, C_{k-1}^r \, C_t^k + {(-1)}^t \; $$
(o\`u  $ \, C_u^v = 0 \, $  si  $ \, v > u \, $)  avec  $ \, r $,
$ t \in \N_+ \, $  et  $ \, r < t \, $:  dans notre cas on a pos\'e  $ t =
n - \big| \Sigma' \big| \, $,  $ \, r = j - \big| \Sigma' \big| \, $  et
$ \, k = d - \big| \Sigma' \big| \, $;  on v\'erifie que l'on a  $ \, j
- \big| \Sigma' \big| < n - \big| \Sigma' \big| \, $  parce que  $ \, j <
n \, $.

\vskip3pt

 (II) \quad  Si  $ \, \Sigma'' \not\subseteq \Sigma' \, $  et
$ \, \big| \Sigma' \cup \Sigma'' \big| > j \, $,  alors l'expression  $
{\big(  E' \big) }_{\Sigma',\Sigma''}$  devient
  $$  {\big( E':2 \big)}_{\Sigma',\Sigma''}
\, = \sum_{\Sb  \Sigma \subseteq \{1, \dots, n\}  \\
\Sigma' \cup \Sigma'' \subseteq \Sigma  \endSb}  \!\! {(-1)}^{n-|\Sigma|}
\,  {(-1)}^{j-|\Sigma'|} \, C_{|\Sigma|-1-|\Sigma'|}^{j-|\Sigma'|} \; .  $$
En regroupant les  $ \Sigma $  qui ont le m\^eme cardinal  $ d $,  un
simple  comptage nous donne
  $$  {\big( E':2 \big) }_{\Sigma',\Sigma''}\, = \sum_{d = |\Sigma' \cup
\Sigma''|}^n
{(-1)}^{n-d} \, {(-1)}^{j-|\Sigma'|} \, C_{d-1-|\Sigma'|}^{j-|\Sigma'|} \,
C_{n - |\Sigma' \cup \Sigma''|}^{d - |\Sigma' \cup \Sigma''|} \; .  $$
   \indent   \`A nouveau, cette derni\`ere expression est nulle gr\^ace
au lemme 3.3, car elle correspond \`a une somme du type  $ \;
\sum\limits_{k=0}^t {(-1)}^{t+r-k} \, C_{k+s}^r \, C_t^k \; $
avec  $ \, r $,  $ t $,  $ s \in \N_+ \, $  et
$ \, r < t \, $. Dans
$  {\big( E':2 \big) }_{\Sigma',\Sigma''}\, $, on a pos\'e
                $ \, t = n - \big|
\Sigma'\cup \Sigma'' \big| \, $,  $ \, r = j - \big| \Sigma' \big| \, $,
$ \, s = \big| \Sigma' \cup \Sigma'' \big| - \big| \Sigma' \big| - 1 \, $
et  $ \, k = d - \big| \Sigma' \cup \Sigma'' \big| \, $;   on v\'erifie que
l'on a  $ \, j - \big| \Sigma' \big| < n - \big| \Sigma' \big| \, $  car
$ \, j < n \, $  et  $ \big| \Sigma' \cup \Sigma'' \big| - \big| \Sigma'
\big| - 1 \geq 0 \, $  car  $ \, \Sigma'' \not\subseteq \Sigma' \, $.

\vskip3pt

 (III) \quad  Si  $ \, \Sigma'' \not\subseteq \Sigma' \, $  et
$ \big| \Sigma' \cup \Sigma'' \big| \leq j \, $,  alors l'expression
$ {\big( E' \big)}_{\Sigma',\Sigma''} $  devient
  $$  {\big( E':3 \big)}_{\Sigma',\Sigma''}
 \, = \sum_{\Sb  \Sigma \subseteq \{1, \dots, n\} \\
\Sigma' \cup \Sigma'' \subseteq \Sigma, \; |\Sigma|>j  \endSb}
{(-1)}^{n-|\Sigma|} \, {(-1)}^{j-|\Sigma'|} \,
C_{|\Sigma|-1-|\Sigma'|}^{j-|\Sigma'|} \; .  $$
Si l'on regroupe les  $ \Sigma $  qui ont le m\^eme cardinal  $ d $,  un
simple  comptage nous donne
  $$ { \big( E':3 \big)}_{\Sigma',\Sigma''}
 = \sum_{d=j+1}^n {(-1)}^{n-d} \, {(-1)}^{j-|\Sigma'|}
\, C_{d-1-|\Sigma'|}^{j-|\Sigma'|} \, C_{n - |\Sigma' \cup \Sigma''|}^{d -
|\Sigma' \cup \Sigma''|} \; .  $$
   \indent   Mais l\`a encore, la derni\`ere expression est nulle d'apr\`es le
lemme 3.3, car elle correspond \`a une somme du type
$$ \; \sum\limits_{k
= j + 1 - |\Sigma' \cup \Sigma''|}^t \!\!  {(-1)}^{t+r-k} \, C_{k+s}^r \,
C_t^k = \sum\limits_{k=0}^t {(-1)}^{t+r-k} \,  C_{k+s}^r \, C_t^k \; $$
o\`u  $ \, C_u^v = 0 \, $  si  $ \, v > u  $, avec  $ \, r $,  $ t $,
$ s \in \N_+ \, $  et  $ \, r < t \, $:  ici on a encore pos\'e  $ \, t = n
- \big| \Sigma' \cup \Sigma'' \big| \, $,  $ \, r = j - \big| \Sigma'
\big| \, $,  $ \, s = \big| \Sigma' \cup \Sigma'' \big| - \big| \Sigma'
\big| - 1 \, $  et  $ \, k = d - \big| \Sigma' \cup \Sigma'' \big| \, $;
toujours pour les  m\^emes raisons,  $ \, j - \big| \Sigma' \big| <
n - \big| \Sigma' \big| \, $   et  $ \big| \Sigma' \cup \Sigma'' \big| -
\big| \Sigma' \big| - 1 \geq 0 \, $.

\vskip3pt

   En conclusion, on a toujours  $ \, {\big( E'
\big)}_{\Sigma',\Sigma''} = 0 $,
d'o\`u  $ \, (E) = 0 \, $,  ce qui termine la preuve.   $ \square $
\enddemo

\vskip23pt

\enddocument